\documentclass[fleqn,10pt]{wlscirep}
\usepackage[utf8]{inputenc}
\usepackage[T1]{fontenc}

\title{Modelling brain-wide neuronal morphology via rooted Cayley trees}

\author[1,2]{Congping Lin}
\author[1]{Yuanfei Huang}
\author[3,4] {Tingwei Quan}
\author[1,2,*]{Yiwei Zhang}
\affil[1]{Center for Mathematical Sciences, Huazhong University of Science and Technology, Wuhan, China}
\affil[2]{Hubei Key Lab of Engineering Modeling and Scientific Computing, Huazhong University of Science and Technology, Wuhan, China}
\affil[3]{Britton Chance Center for Biomedical Photonics, Wuhan National Laboratory for Optoelectronics-Huazhong University of Science and Technology, Wuhan, China}
\affil[4]{MoE Key Laboratory for Biomedical Photonics, Collaborative Innovation Center for Biomedical Engineering, School of Engineering Sciences, Huazhong University of Science and Technology, Wuhan, China}

\affil[*]{Correspondence to yiweizhang@hust.edu.cn}


\begin{abstract}
Neuronal morphology is an essential element for brain activity and function. We take advantage of current availability of brain-wide neuron digital reconstructions of the Pyramidal cells from a mouse brain, and analyze several emergent features of brain-wide neuronal morphology. We observe that axonal trees are self-affine while dendritic trees are self-similar. We also show that tree size appear to be random, independent of the number of dendrites within single neurons. Moreover, we consider inhomogeneous branching model which stochastically generates rooted 3-Cayley trees for the brain-wide neuron topology. Based on estimated order-dependent branching probability from actual axonal and dendritic trees, our inhomogeneous model quantitatively captures a number of topological features including size and shape of both axons and dendrites. This sheds lights on a universal mechanism behind the topological formation of brain-wide axonal and dendritic trees.

\end{abstract}
\begin{document}

\flushbottom
\maketitle

\thispagestyle{empty}

\section*{Introduction}
Neurons, the primary components of central nervous system are electrically excitable cells that receive, process, and transmit information through electrical and chemical signals between each other. Digital reconstructions of neurons provide information for quantitative measurements of neuronal morphology. With digital reconstructions, many studies have highlighted the importance of neuronal morphology in its biological function. Mainen and Sejnowski \cite{Mainen1996} have shown a causal relationship between dendritic structures and intrinsic firing patterns observed from {\em in vitro} electrical recordings for a wide variety of cell types. Vetter {\em et al.} have shown that branching pattern strongly affected the propagation of action potentials which links information processing at different regions of the dendritic tree \cite{Vetter2001}. Ferrante {\em et al.} have shown that even subtle membrane readjustments at branch point could drastically alter the ability of synaptic input to generate and propagate the action potentials \cite{Ferrante2013}. It has also been shown that the input-output response function of neuron’s dendritic arbor grows with the tree size \cite{Gollo2013,Publio2012}. More recently, Yi {\em et al.} have demonstrated a crucial role of neuronal morphology in determining field-induced neural response \cite{Yi2017}. These results, among others, have contributed to a now widespread acceptance that neuronal morphology plays a critical role in its activity and function.

The formation of neurons through branching is driven by complex interactions of intracellular and extracellular signaling cascades which are proving difficult to be completely understood by molecular biology alone. Mathematical or computational modelling approaches instead provides an alternative and complementary approach to uncover mechanism underlying neuronal morphology. Simulators {\em L-Neuron} \cite{Ascoli2000,Ascoli2001} and NeuGen \cite{Eberhard2006} were developed to create virtual neuronal structures in silico by means of iteratively sampling experimental statistical distributions of shape parameters (including e.g. branch diameter, length, ect). In contrast to L-Neuron or NeuGen simulators which require a large number of experimental neuron samples to obtain reliable distributions of shape parameters, Van Ooyen {\em et al.} used models in which the morphology of a single dendritic tree was represented in a highly abstract manner \cite{VanPelt1985,VanPelt1986a,VanPelt2004,vanPelt2007} where growing neurons were modeled as expanding, circular neuritic fields. Based on such phenomenological dendritic growth models, NETMORPH was developed to simulate 3D neuronal networks from the perspective of individual growth cone \cite{Koene2009}, using simple rules for neurite branching probability at each terminal segment. Other modelling approach with simple rules to account for the spatial embedding of tree structures, includes optimal rewiring and particle-based diffusion limited aggregation approach: the optimal rewiring approach generates branching geometry by minimizing the wiring length and the path-length to root between branch points and synapses in dendrite trees \cite{Cuntz2010,Cuntz2012}; the particle-based diffusion limited aggregation approach provides model-based measures to estimate ``diffusive'' shape of neuronal tree-like structures \cite{Luczak2010, Luczak2006}.

However, analysis and modelling of neuron morphology to the best of our knowledge so far are largely based on neuron digital reconstructions from certain regions or layers of a brain (e.g. \cite{Ascoli2000,Koene2009}). Moreover, a majority of these neuronal morphology analysis and modelling concentrate on neuron dendrites, with only a few studies of axonal branching structure in certain layers recently (e.g. \cite{Gillette2015b, Mohan2015}). Knowledge about how dendrites and axons branch in an entire brain is still limited. Recently, breakthroughs in imaging \cite{Li2010,Ragan2012} and molecular labeling \cite{Jefferis2012,Chung2013} techniques have provided tools to trace and digital reconstruct the almost complete morphology of neuronal populations at a single-axon resolution through a whole brain \cite{Zhou2018,Quan2016}; this offers opportunity to quantify brain-wide dendritic and axonal branching morphology. In this manuscript, we take advantage of current availability of brain-wide neuron digital reconstructions of pyramidal cells from a mouse brain \cite{Zhou2018} and analyze the complete axonal and dendritic branching morphology. In particular we show that axons are topologically self-affine whereas dendrites are topologically self-similar, and topological structures of both axon and dendrites are far away from symmetry and appear to be random. We also develop a self-organized probabilistic model for the entire axonal and dendritic branching structures. In contrast to growth models \cite{vanPelt2007} which require a simulation/growth time to create virtual trees of finite size, we use the generation of rooted 3-Cayley trees (where non-terminal nodes are linked to 3 neighbors; in neuron trees, 3 linkages of a node represent one mother branch and two daughter branches) via a stochastic branching process. Using estimated branching probability from brain-wide axonal and dendritic trees, we show that this simple probabilistic model is sufficient to quantitative recapture several statistical properties of neuronal morphology including distributions of neuron size and topological width/length.

\begin{figure}[h!]
	\centering
	\includegraphics[width=18cm]{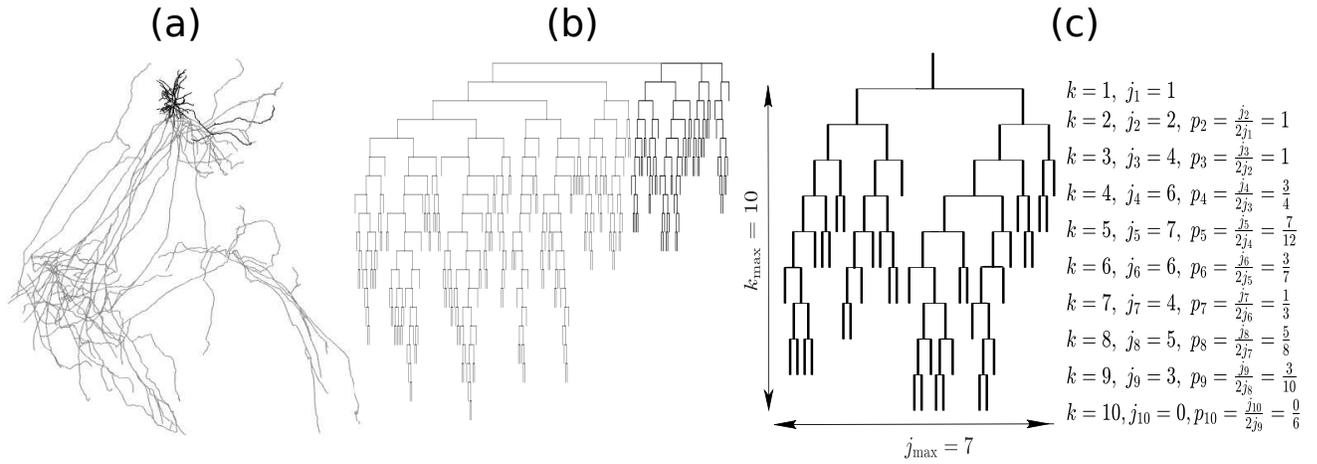}
	\caption{{\bf Illustration of one single brain-wide neuron digital reconstruction}. (a) shows an example of neuronal digital reconstruction; data is from \cite{Zhou2018}. (b) shows the corresponding Dendrogram of the neuron shown in (a) based on its the topological branching structure. The axon part is colored in gray while dendrite part colored in black in both panels (a) and (b). (c) shows the dendrogram based on one dendrite. The number of nodes $j_k$ and branching frequency $p_k$ at each order $k$ is illustrated in this example. $k_{\max}$ and $j_{\max}$ are the maximum values of $k's$ and $j's$ and are referred as tree length and width respectively. In this example, the maximal order equals to $10$ (i.e, $k_{\max}=10$) and maximal number of nodes at different orders is $7$ which appears at order $k=5$ (i.e., $j_{\max}:=\max_k\{j_k\}=7$). Note that the total number of nodes at order $k$ equals to $2j_{k-1}$.}
	\label{fig_eg}
\end{figure}

\section*{Results}
\subsection*{Brain-wide axonal and dendritic morphology quantifications}
In this section, we analyse neuron morphology from 35 brain-wide digital reconstructions of pyramidal neurons in a mouse brain \cite{Zhou2018}. A typical neuron consists of a cell body (soma), dendrites, and an axon. Here we analyze axon and dendrite tree morphology in aspects of size, asymmetry, and shape of trees as well as correlation among axon and dendrites in single neurons.

\subsubsection*{Topological structure of axonal and dendritic arbors}
We first focus on topological structure of dendritic and axonal branching from soma to neuron terminals, despite the spatial position of each branch. Fig~\ref{fig_eg} illustrates an example of a single brain-wide neuronal morphology and its topological branching structure. We seperate individual neurons into axon and dendrites for analysis; see Fig~\ref{fig_eg} as an example. From 35 actual neuron reconstructions, we find that on average there are $6.8\pm 0.56$ number of dendrites besides one axon. Neurons have been widely viewed as binary trees (e.g.\cite{Vormberg2017}), here we first examine this by looking into the degree (i.e. the number of branches each node links to) of each node. Clearly, each terminal is connected to one branching node (i.e. degree 1) and branching nodes are mainly linked with three branches in both axons and dendrites - with only $<1\%$ of branching nodes in axons and $2\%$ in dendrites linking to more than 3 branches; those could be due to limited spatial resolution in imaging which leads to possible errors in digital reconstructions. Overall, this suggests that brain-wide neurons are well approximated as binary trees.

\begin{figure}[h!]
	\centering
	\includegraphics[width=16cm]{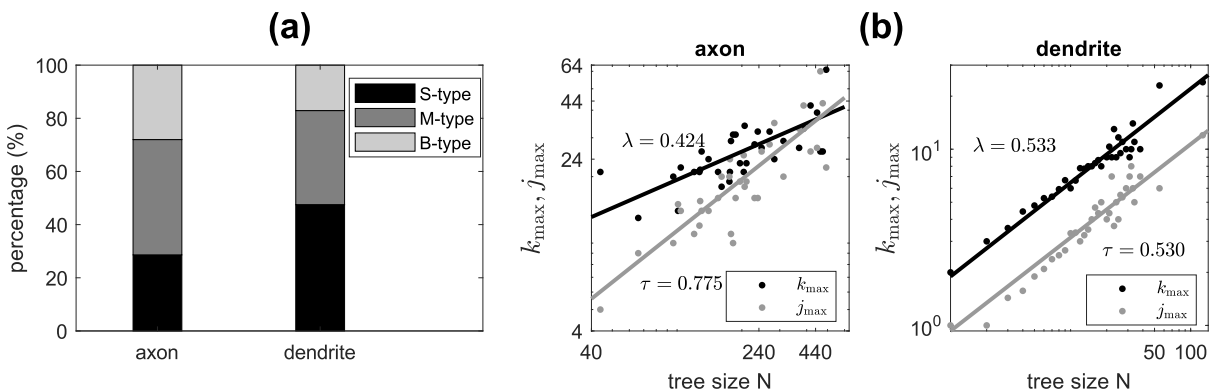}
	\caption{{\bf Tree branching nodes and topological shape.} (a) shows proportion of each type of branching nodes in axons ($n=35$) and dendrites ($n=239$) from 35 actual brain-wide neuron reconstructions. (b) shows topological length and width and their scaling with tree size as $k_{\max}\sim N^{\lambda}$ and $j_{\max}\sim N^{\tau}$; each dot represent the $k_{\max}$ (black) or $j_{\max}$(grey) averaged over trees with the same tree size $N$. Exponents $\lambda$ and $\tau$ are indicated in each panel for both axons and dendrites. They show significant difference in axon (F-test; $p=0.012$) and no significant difference in dendrite (F-test; $p=0.97$).}
	\label{fig_scaling_exp}
\end{figure}

We quantify the number of branching nodes (which is equivalent to the number of terminals minus one) for both axonal and dendritic trees, and refer as the topological size. On average, we find that axons are of size $224.1 \pm  20.89 (n=35)$, whereas dendrites are significantly smaller, of size $8.29\pm 0.72 (n=239)$, i.e. about $1/27$ of the axonal tree size. We also classify the branching nodes into three types~\cite{Gillette2015a}: B-type if both child branches themselves bifurcate; M-type if only one child branch bifurcate; S-type if branching with two terminals. By such classification, we find that axons have smaller proportion of S-type nodes whereas dendrite trees have larger proportion of S-type branching nodes; see Fig~\ref{fig_scaling_exp}~(a).

We characterize the topological shape of neuron trees by measuring its topological width and length as illustrated in Fig~\ref{fig_eg}~(c). To do this, we first give an order (denoted as $k$) for each node as its topological distance to the soma, i.e. the number of branches in its path to the soma. The topological length (denoted as $k_{\max}$) of a tree is then defined as the maximal order among nodes. For the topological width, we calculate the number of branching nodes (denoted as $j_k$) for each order $k$ and define the topological width (denoted as $j_{\max}$) of a tree as the maximum among all $j_k, k=1,\cdots,k_{\max}$ \cite{Molnar2006}. This is similar to Sholl analysis \cite{Sholl1953} when topological distance is applied. In the population of our 35 actual brain-wide neuron reconstructions, we observe that on average $j_{\max}\propto N^{\tau}, k_{\max}\propto N^{\lambda}$ and the scaling exponents are different between axon ($\tau=0.775 \pm 0.114, \lambda=0.424\pm 0.079$) and dendrite ($\tau=0.530 \pm 0.035, \lambda=0.533\pm 0.036$); see Fig~\ref{fig_scaling_exp}·(b). Note that in axons, exponents $\tau$ and $\lambda$ are significantly different (F test; $p=0.012$) where in dendrites there is no significant difference between the exponents (F-test $p=0.97$); this  suggests that in aspects of topological shape, dendrites are approximately self-similar (i.e. $\tau\sim\lambda$) whereas axons are self-affine (i.e. $\tau \nsim \lambda$).

Next, we quantify morphological asymmetry of neuron trees. Following the definition in \cite{VanPelt1992}, tree asymmetry $A$ describes the average of local partition asymmetry $A_p(r_j,s_j):=\frac{|r_j-s_j|}{r_j+s_j-2}$ for all branching nodes (here $r_j,s_j$ are the number of terminals in two subtrees from a branching node) except for those with partition $(1,1)$ (i.e. branching to two terminals). This tree asymmetry reflects the relative balance of branching within a tree based on the distribution of terminals between two subtrees birfucated from a node. This asymmetry index ranges from 0 (symmetry) to 1 (asymmetry). The axon and dendrite from brain-wide neuron digital reconstructions show no significant difference (student t-test $p=0.274$) on tree asymmetry with $A=0.791 \pm 0.006 (n=35)$ for axonal trees and $A=0.743 \pm 0.019 (n=211)$ for dendritic trees and both are far away to symmetry. These measured partition asymmetry values are close to reported values in \url{http://neuromorpho.org/}. Another asymmetry measurement - excess partition asymmetry introduced by Samsonovich and Ascoli \cite{Samsonovich2006} considers the difference between partition asymmetry actually measured at one branch and the average of partition asymmetry computed for the same branch after all possible shuffling of the granddaughter branches. This excess partition asymmetry is useful to test the randomness in branching or the existence of control/regulation process within neurons. Interestingly, in contrast to 'regulation mechanism' suggested in \cite{Samsonovich2006} based on systematically positive excess partition asymmetry $E_p$ measured from pyramidal cells on some slice of a brain \cite{Samsonovich2005}, we show that our brain-wide axonal and dendritic trees appear to have almost no excess partition asymmetry with $E_p=0.012 \pm 0.005 (n=35)$ (t-test to zero mean with $p=0.013$) for axon and $E_p=0.004 \pm 0.006 (n=151)$ (t-test $p=0.54$). This suggests branching in both axonal and dendritic trees appears to be random.

\subsubsection*{Geometric size of axon and dendrites}
\begin{figure}[h!]
	\centering
	\includegraphics[width=12cm]{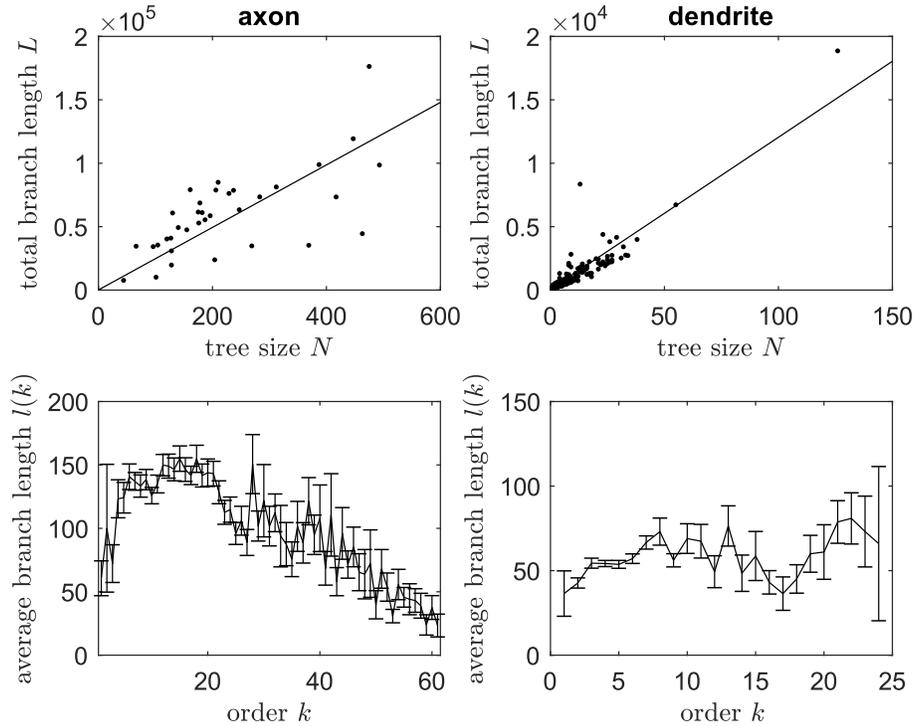}
	\caption{{\bf Tree geometric length.} Top two panels show total branch length $L$ (measured as the sum of all branch length) as a function of topological tree size $N$; lines give the best fitting to the function $L(N)=m(2N+1)$ with $m=123.2\mu m$ for axons and $m=59.9\mu m$ for dendrites. Bottom two panels show the average branch length (mean $\pm$ sem is shown) at each order $k$ for axons (left) and dendrites (right).}
	\label{fig_length}
\end{figure}

Regarding geometric size of neuron trees, we measure the total branch length $L$ of a tree (i.e. the sum of all branch length in a tree). On average, axons have a total branch length of $5.97\times 10^4\pm 0.55\times 10^4\mu m$ (median $5.86\times 10^4 \mu m$) and dendrites have a total branch length $9.6\times 10^2\pm 0.99\times 10^2 \mu m$ (median $5.6\times 10^2\mu m$). Moreover, we observe that the total length $L$ strongly correlates with tree size $N$ (Pearson $r=0.68, p<0.0001$ for axons and $r=0.91, p<0.0001$ for dendrites) and increase approximately linear as the tree size $N$ for both axons and dendrites; see Fig~\ref{fig_length}~top panels. The mean total length is expected to be the mean branch length (say $m$) multiplied by the number of branches $2N+1$ for a given tree size $N$. The data show a best linear fit to the function $L=m(2N+1)$ with a mean branch length $m=123.2\mu m$ for axons and $m=59.9\mu m$ for dendrites. Such strong correlation and similar linear slope between axon length and its tree size were also observed in vitro axons of human neurons across all layers of medial temporal cortex, though total length of axons in human brains is much larger than that in mouse \cite{Mohan2015}; this suggests a potential universal mechanism controlling branch length cross species. We also note that in our brain-wide neurons, though the total branch length approximately grows linearly with the number of branches, the branch length is not uniform along orders; see Fig~\ref{fig_length}~bottom panels for both axons and dendrites. In particular, the branch length of axons is larger at intermediate orders than those at smaller or larger branch orders.

\subsubsection*{Correlation among axon and dendrites in single neurons}
\begin{figure}[h!]
	\centering
	\includegraphics[width=16cm]{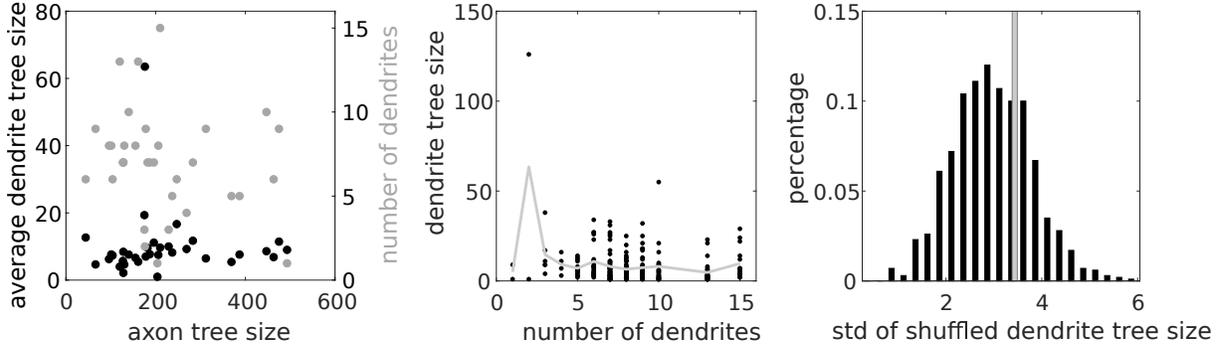}
	\caption{{\bf Tree size of axons and dendrites from single neurons.} (a) shows scatter plots of average dendrite tree size and number of dendrites against axon tree size from single neurons; each dot represents measurements from one neuron. Pearson correlation analysis suggests no significant correlation with $p=0.98$ for average dendrites tree size and $p=0.23$ for number of dendrites. (b) shows the tree size of dendrites in relation with the number of dendrites in single neurons; each dot represents measurements from one dendrite. Gray lines represent the average cross vertical points for each number of dendrites. Pearson correlation analysis suggests no significant correlation between the number of dendrites and average dendrite size ($p=0.061$). (c) shows the histogram of standard deviations of average dendrite tree sizes from 1000 random shuffling of dendrites among 7 neurons (each of which has 7 dendrites); gray bar shows the standard deviation of average dendrite tree sizes from the actual 7 neurons.}\label{fig_random}
\end{figure}

It has been suggested that in one given portion of a neuron, dendritic morphology may be under intrinsic homeostatic control which regulates tree size fluctuations systematically by counterbalancing the remaining dendrites in the same cell \cite{Samsonovich2006}. We test such a control mechanism in axons and dendrites from our actual brain-wide neurons. If trees of a neuron were mainly up(down)-regulated by common factors, then neurons with larger axons in terms of tree size would have larger dendrites; if axons and dendrites were mutually regulated by competition, then neurons with larger axons would have smaller dendrites. Based on calculated tree size of axons, average tree size among dendrites and the number of dendrites for each of 36 neurons, we perform Pearson correlation analysis among these quantities. Interestingly, we find no significant correlation between axon tree size and average dendrite tree size ($p=0.98$), between axon tree size and number of dendrites ($p=0.23$), or between average dendrite tree size and the number of trees ($p=0.061$). These suggest that axon tree size, dendrite tree size and dendrite number in single neurons are likely to be independent. Moreover, if tree size is regulated within single neurons, then tree shuffling would lead to an increase in the variance of average tree size in the neuron population; if tree size is random among populations, then shuffling would not lead to significant change on the variance \cite{Samsonovich2006}. To further test such randomness, we consider neurons with the same number of dendrites. Then we randomly shuffle the dendrites in these neurons and calculate the standard deviation (std) of the average dendrite tree size among single neurons. Fig~\ref{fig_random}~(c) show that 35.5\% of 1000 random shuffling gives higher std of the average dendrite tree size than the std from actual neurons, suggesting tree size of dendrites within single neurons is random.


\subsection*{Modelling topology of axonal and dendritic trees}
As indicated by above statistical analysis, the tree size of dendrites and axons appear to be random independent of the number of dendrites, indicating that the formation of brain-wide neuron trees is likely to be a stochastic process. In this section, we use stochastic generation of rooted z-Cayley trees \cite{Ostilli2012} to model the topological structure of axonal and dendritic arbors. A z-Cayley tree is a tree where non-terminals have z-neighbours. Here to model neuron trees, we use $z=3$ as non-terminal nodes in neuron trees are shown to mainly have 3 branches (3 neighboring nodes): one in previous order and two in the next order. In the generation of a rooted 3-Cayley tree, we start with an arbitrary node as the root (which represents the soma in a neuron tree), and assign an link from this root connecting to another node (viewed as the 1-th order node); the generation of tree proceeds from this 1-th order node, by 1) bifurcating to two daughter branches (linking to 2 neighbouring nodes in the next order) and 2) randomly selecting its neighboring nodes in the next order with certain probability; if one node in the $k-$th order is selected, then it bifurcate into two daughter branches to next order $k+1$ and random selection on its two neighbours applies independently; if a neighbouring node is not selected, then this neighboring node is considered as a terminal (which could not bifurcate into the next order). The procedure continues until no more nodes are selected to bifurcate and eventually generates a rooted 3-Cayley tree.

We first investigate homogeneous branching model where all nodes are selected with the equal probability $p$. In this case it is well known that beyond a critical probability $p_c=1/2$, a positive probability exists that an infinite tree is generated, whereas below this threshold, the probability for such infinite spanning is zero and a finite tree is generated \cite{Albert2002}. Moreover, the mean tree size of generated finite trees increases with branching probability as
\begin{equation}\label{eq_N-p}
\langle N\rangle=1/(1-2p),~\mbox{~~if~~} p<p_c=1/2.
\end{equation}
In fact if we denote $T$ by the mean number of branches in one branch from the 1-th order node, then by statistical equivalence between sub-branches $T=p(1+2T)$, which gives $T=p/(1-2p)$, and the mean tree size $\langle N \rangle=1+2T$ where $1$ corresponds to the 1-th order branching node.

\begin{figure}[h!]
	\centering
	\includegraphics[width=14cm]{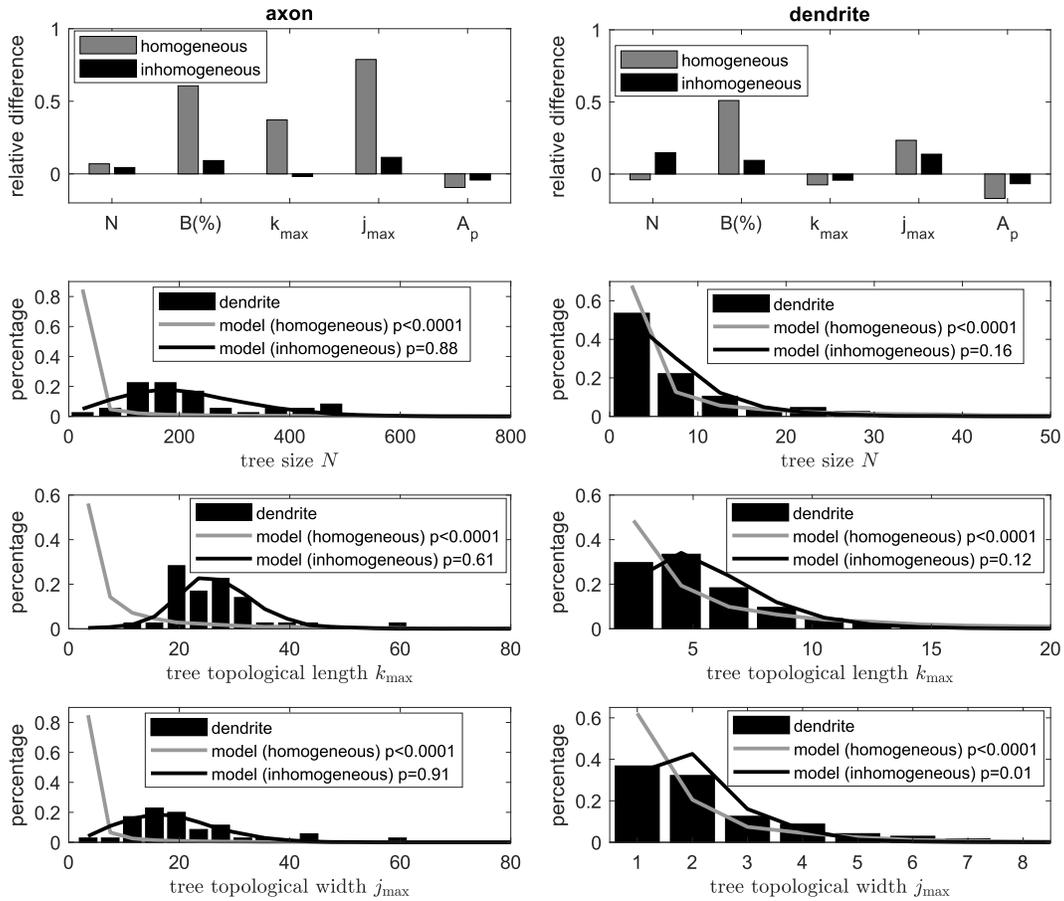}
	\caption{{\bf Comparison between actual neuron data and model.} Top two panels show relative difference on quantities including tree size $N$, percentage of B-type branching nodes $B(\%)$, tree topological length $k_{\max}$ and width $j_{\max}$ as well as tree asymmetry $A$ from (in)homogeneous modelling to those measured from actual axonal (left) and dendritic (right) trees respectively. Bottom panels show  distributions on tree size $N$ and topological shape parameters ($k_{\max}$ and $j_{\max}$). Two-sample Kolmogorov–Smirnov test is used to compare between distributions from actual axonal (left) and dendritic (right) trees and virtual trees generated using the (in)homogeneous model; $p$ values are indicated for each comparison. The imhomogeneous model uses the best fitted exponential decay with non-zero plateau of the branching probability given in Fig~\ref{fig_exp_branching}.}
	\label{fig_model_data}
\end{figure}

To model the axon and dendrite trees, we estimate the constant branching probability $p$ from Eq.~\eqref{eq_N-p} using measured average tree size, which gives $p=0.498$ for axon and $p=0.44$ for dendrite. Indeed, the model with estimated branching probability produce similar mean tree size as the actual data; see Fig~\ref{fig_model_data}~top two panels. We then test this simple homogeneous model for other quantities described in the above section, such as proportions of 3 types of branching nodes (we take the B-type branching node as an example for the comparison), topological shape indexes ($j_{\max}$ and $k_{\max}$), and tree asymmetry $A$. Note that the model takes independent probability for each node and thus no excess asymmetry is expected. Fig~\ref{fig_model_data}~top panels show the relative difference of the modelling results to actual data on these quantities; in particular, we see that proportion of B-type branching node and topological width ($j_{\max}$) from homogeneous model from both axons and dendrites are far different ($>20\%$ relative difference) to measurements from actual data. Moreover, distributions on tree size ($N$), tree topological length ($k_{\max}$) and width ($j_{\max}$) as shown in Fig~\ref{fig_model_data}~bottom panels for both axon and dendrite are significantly different to corresponding distributions from actual axonal and dendritic trees ($p<0.0001$ for both axons and dendrites). We also notice that distributions from the homogeneous model illustrate higher frequency at lower values and larger variation on these quantities for both axon and dendrite, indicating that it has high probability to generate smaller tree than the expected mean. We thus conclude that homogeneous model is not sufficient to model the topology of neuron trees.

\begin{figure}[h!]
	\centering
	\includegraphics[width=8cm]{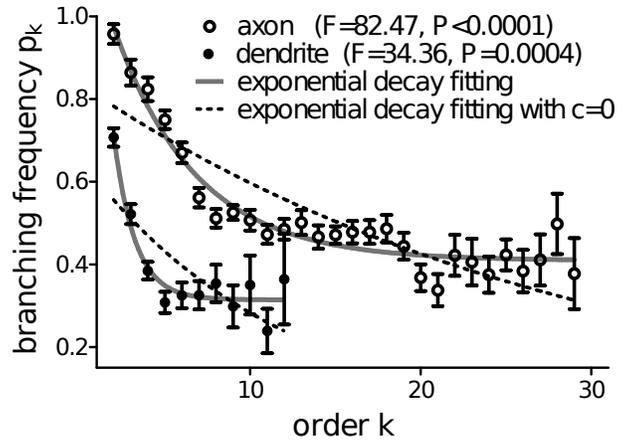}
	\caption{{\bf Branching frequency of axonal and dendritic trees.} F test shows that the branching frequency measured from actual axonal (circles) and dendritic (dots) trees prefers an exponential decay $b\exp(-ak)+c$ with non-zero plateau $c$ (solid lines), instead of an exponential decay to $c=0$ (dashed lines). The best fitting parameters for the exponential decay with non-zero plateau are $a=0.206, b=0.855, c=0.409$ for axonal trees and $a=0.79, b=1.933, c=0.313$ for dendritic trees. To have appropriate average branching frequency values for the fitting, only orders with sample size larger than 10 is considered here.}
	\label{fig_exp_branching}
\end{figure}

We next consider inhomogeneous model where branching probability differs between nodes. In particular we investigate order dependent branching probability $p_k$; such dependence was considered in the growth model\cite{VanPelt2004} and the computational model \cite{Donohue2008}. Fig~\ref{fig_exp_branching} indeed shows that the branching probability does depend on the order in particular at low orders. In the growth model \cite{VanPelt2004} for dendritic topology, branching probability is assumed to be decay exponentially to $0$ as the  order $k\to \infty$. However, the observed branching frequencies for both axon and dendrites in our brain-wide neurons do not fit well to exponential decay with zero plateau, instead, the data exhibit steady branching frequencies at high orders $k$ and statistically the data fit better to exponential decay with non-zero plateau in the form of $p_k=b\exp(-ak)+c$ (F-test; $p<0.001$ for both axons and dendrites). In the best fitting curve, as $k\to \infty$, the branching probability $p_k\to c<p_c=0.5$ (the critical probability); this ensures a finite tree generation in our model \cite{Albert2002}. The plateau $c$ also reflects the almost constant branching probability at large orders where subtrees branch in an approximately ``homogeneous'' manner.

\begin{figure}[h!]
	\centering
	\includegraphics[width=14cm]{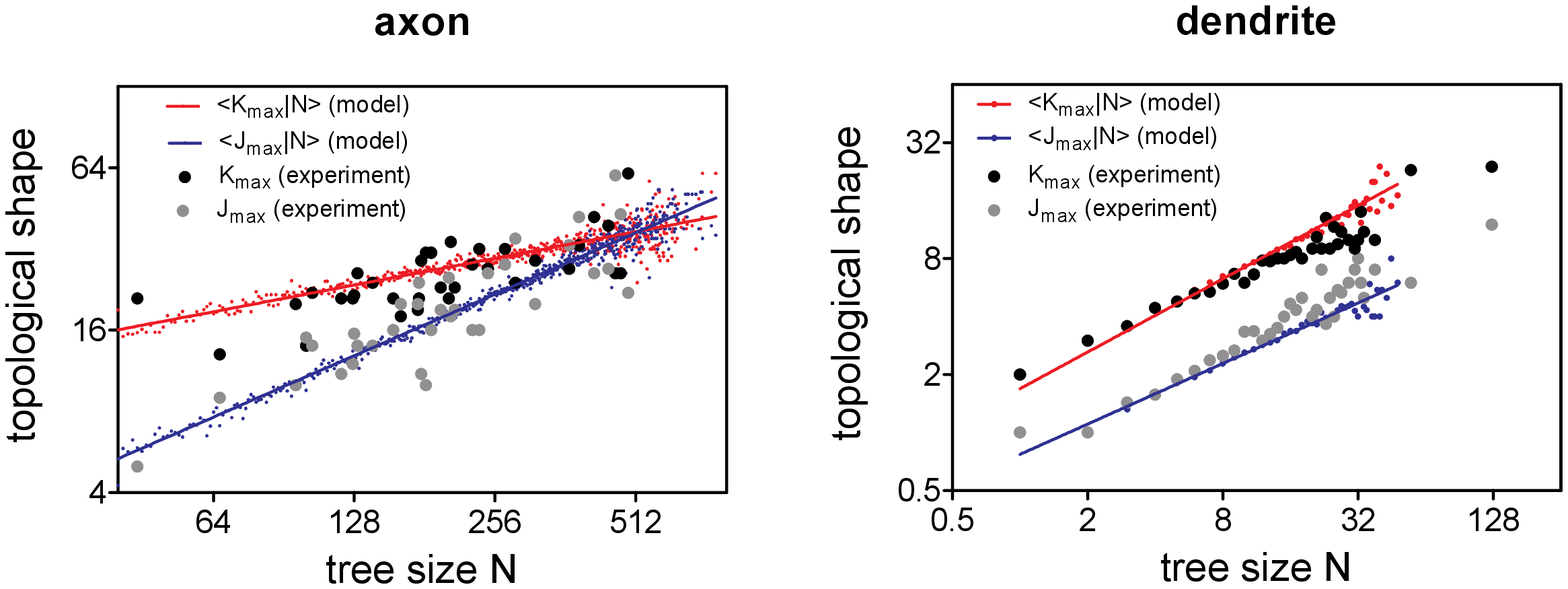}
	\caption{{\bf Tree topological width and length scale with tree size from model simulation.} Red and blue dots show the conditional average of topological length $k_{\max}$ and width $j_{\max}$ respectively for given tree size (i.e. $\langle k_{\max}|N\rangle$ and $\langle j_{\max}|N\rangle$) from simulations using inhomogeneous model with parameters given in Fig~\ref{fig_exp_branching} for both axonal (left panel) and dendritic (right panel) trees. Lines show the best fitted curve of the corresponding simulated dots to $\langle k_{\max}|N\rangle\sim N^{\lambda}$ and $\langle j_{\max}|N\rangle\sim N^{\tau}$. Best fitting gives exponents $\lambda=0.339\pm 0.007, \tau=0.754\pm 0.009$ (significant difference between $\lambda$ and $\tau$; F-test; $p<0.0001$) for axons and $\lambda=0.631\pm 0.048, \tau=0.522\pm 0.044$ (no significant difference between $\lambda$ and $\tau$; F-test; $p=0.339$) for dendrites. Black and grey dots are measured from actual data as shown in Fig~\ref{fig_scaling_exp}~(b).}
	\label{fig_scaling_model}
\end{figure}

With the best fitting branching probability in the form of exponential decay to a non-zero plateau, the inhomogeneous model leads to mean tree sizes $\langle N\rangle\approx 214.3$ and $7.1$ for axonal and dendritic trees respectively, which are close to measured values from actual data. We also test this inhomogeneous model on other topological quantities and show in Fig~\ref{fig_model_data} that the model not only recaptures the mean tree topological measurements (with relative difference $<20\%$ to measurement from actual neurons), but also agrees well on their distributions in particular the tree size and topological shapes $k_{\max}$ and $j_{\max}$ (two-sample Kolmogorov–Smirnov test with $p\geq 0.05$, except for test on $j_{\max}$ of dendrite where $p=0.01$). We next examine the scaling of topological shapes in relation to tree size. Fig~\ref{fig_scaling_model} shows that the conditional averages of tree length and width scale with the tree size $N$ as $\langle k_{\max}|N\rangle\sim N^{\lambda}$ and $\langle j_{\max}|N\rangle\sim N^{\tau}$ with exponents $\lambda=0.339\pm 0.007, \tau=0.754\pm 0.009$ (significant difference between $\lambda$ and $\tau$; F-test; $p<0.0001$) for axon and $\lambda=0.631\pm 0.048, \tau=0.522\pm 0.044$ (no significant difference between $\lambda$ and $\tau$; F-test; $p=0.339$) for dendrite. These exponents from simulations are close to those estimated from actual neurons, showing that as similar as in actual neurons, simulated virutal axons are self-affine while simulated virtual dendritic trees are self-similar. These suggest that inhomogeneous model with a simple order dependent rule on branching probability is sufficient to statistically capture a number of topological features of both axonal and dendritic morphology in a quantitative manner.

\subsection*{Statistical properties of the probability model}

In this section, we explore how topological characteristics of artificial trees generated from the inhomogeneous model with branching probability $p_k=\min\{b\exp(-ak)+c,1\}, \forall k\geq 2$ vary with parameters. Note that we take a $\min$ operation in $p_k$ here, considering a probability cannot be larger than $1$. Also note that, if we denote $p_1$ as the branching probability for order $k=1$, then $p_1\equiv 1$, as the first branching is initially selected in the model.

\begin{figure}[h!]
	\centering
	\includegraphics[width=18cm]{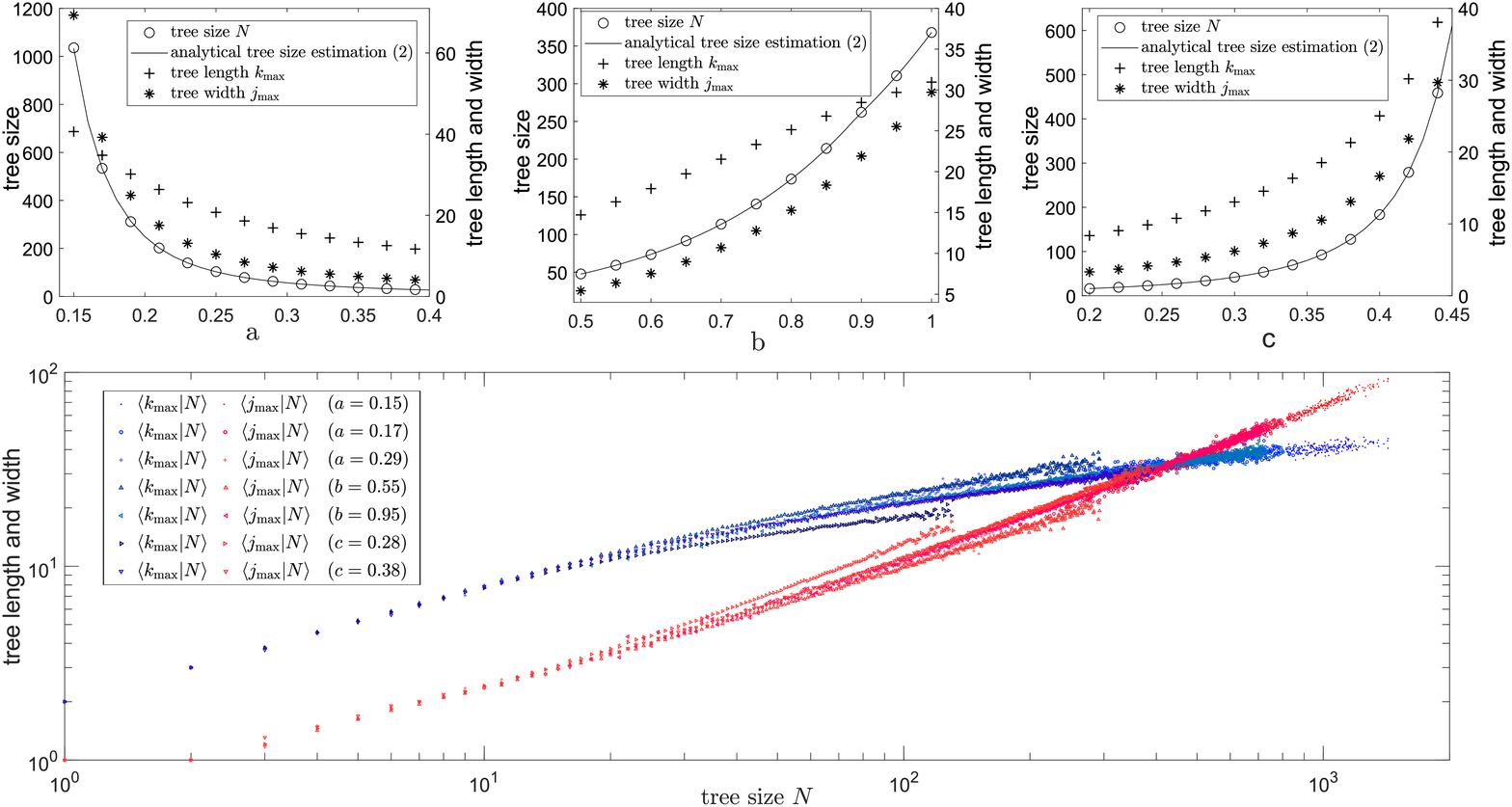}
	\caption{{\bf Statistical properties of tree topological shape when varying model parameters.} Top panels show average tree size, tree length and width when varying parameters $a, b, c$ in the inhomogeneous model with branching probability $p_k=\min\{b\exp(-ak)+c,1\}$ from simulations; lines show the mean tree size according to estimation~\eqref{eq_treesize}. Bottom panel shows the change of conditional average of tree length and width $\langle k_{\max}|N\rangle, \langle k_{\max}|N\rangle$ using different parameters in the inhomogeneous model; one parameter is given in the legend while other parameters are the same as the best fitting curve for axons in Fig~\ref{fig_scaling_exp}~(b); for a given tree size $N$, only sample size larger than 10 of tree length and width are plotted.}
	\label{fig_scaling_model_abc}
\end{figure}

Simply, the mean tree size which can be calculated as the total number of branching nodes from all orders, and read as
\begin{equation}\label{eq_treesize}
\langle N\rangle=\sum_{k=1}^{\infty}j_{k}=1+\sum_{k=2}^{\infty}2^{k-1}\prod_{m=2}^{k}p_m
\end{equation}
For homogeneous model $p_m=p$, it reduces to Eq~\eqref{eq_N-p}. For inhomogeneous model, this estimation agrees well with simulated data; see Fig~\ref{fig_scaling_model_abc}~top panels. As expected, tree size $N$, tree length $k_{\max}$ and width $j_{\max}$ increase with decreasing decay rate $a$ or increasing $b, c$. Moreover, we note that for a similar variation of parameters, the decay rate $a$ allows to give a board range of values on topological quantities in particular the tree size $N$. Furthermore, from the simulated data, we observe that on average tree length $k_{\max}$ is larger than tree width $j_{\max}$ when tree size is small, whereas when tree size is big (e.g. see simulations with small $a$ on the top left panel in Fig~\ref{fig_scaling_model_abc}) tree length $k_{\max}$ is smaller than tree width $j_{\max}$ on average. This is further confirmed in Fig~\ref{fig_scaling_model_abc}~bottom panel, which illustrates changes of conditional average of tree length $k_{\max}$ and width $j_{\max}$ for various tree size $N$. Note that the model has three parameters $a,b$ and $c$ and different combinations of the parameters could give a same tree size $N$ depending on the branching probability express on $p_k$. However, from Fig~\ref{fig_scaling_model_abc}~bottom panel, we see that for the same tree size $N$, the corresponding mean tree length and width from different combinations of parameters are similar. This suggests that the tree shape largely depends on the tree size regardless of the branching probability pattern along the orders in the model. Moreover, Fig~\ref{fig_scaling_model_abc}~bottom panel shows that for small magnitude of tree size $N$, the mean tree length and width grow with the tree size $N$ in a similar scale, i.e. $\langle k_{\max}|N\rangle\sim N^{\lambda}, \langle j_{\max}|N\rangle\sim N^{\tau}$ with $\lambda\approx \tau$; the dendrite trees are particular of this case. Whereas for a large magnitude of tree size $N$, the mean tree width grow faster than tree length when increasing tree size $N$, i.e. $\lambda<\tau$; axons are particular of this case. An intersection between tree length and tree width in relation with tree size occurs at around $N=400$; for tree size much larger than this, the mean tree width is larger than mean tree length (i.e. $\langle j_{\max}|N\rangle>\langle k_{\max}|N\rangle$). These could be used as predictions for topological tree length and width when knowing tree size.

\section*{Conclusion and Discussion}
In this manuscript we examine several neuronal morphological features and develop an inhomogeneous model to generate finite realistic virtual axonal and dendritic tree structures. By analyzing topological shape (length and width) in relation with tree size, we find that axons exhibit an self-affine patten while dendrites are self-similar. Moreover, analyses on excess asymmetry and dendrite shuffling, suggest that tree size appear to be random; this supports the idea of a probability model for neutron tree structures. Using inhomogeneous branching probability, our probability model generates finite virtual trees of statistically similar to brain-wide Pyramidal neutrons in a number of features (including axonal and dendritic tree size and shape). In contrast to published simulation tools \cite{Eberhard2006,Ascoli2000,Ascoli2001,vanPelt2007}, this modelling approach provides a tool with simple rules in a self-organized manner to generate virtual topological neuron structures of finite size. We remark here that models using different order dependent branching probability functions (e.g. polynomial decay, modelling data not shown) which fit the measured branching frequency as shown in Fig~\ref{fig_exp_branching}, are also able to recapture measured data. Mechanisms underlying the decay pattern of branching probability remains unclear; however order dependent branching frequency alone is thus sufficient to capture the tree topology.

Our modelling of brain-wide neuron branching topology can be extended to take geometric information of branches into account, e.g. branch length and branch angles as studied in \cite{Koene2009,Donohue2008}. Note that a branch in a real neuron is not a straight line between two branching points, but exhibit complicated curvature structures; see Fig~\ref{fig_eg}~(a). Including curvature structures would thus be important for a complete and better understanding of mechanism underlying brain-wide neuron morphology and its relation to specific functions at a single neuron scale. One possibility for such extending could be the inclusion of a stochastic process for the outgrowth direction of each branch as considered in the NETMORPH \cite{Koene2009}. Single neurons can expand towards different regions in the brain and different areas in the brain are specialized for different functions; e.g. left hemisphere of the brain is dedicated to language while the right hemisphere is involved in more creative activities such as drawing. Thus morphological difference between different regions in a single neuron scale is likely to link at with different region-specified functions \cite{Zeng2017,Joo2013,Han2018}. Indeed, neuron morphology and function in different layers of brains have been studied \cite{Schubert2006,Marx2013}. However, we are still far away to understand neuron morphology-function relationship at a single neuron scale.

As a theoretical approach in studying morphological role of neurons in their function, our brain-wide modelling of neuron branching could be useful to  study the enhancement of dynamic range \cite{Gollo2009}, by creating virtual realistic neuron topologies. Traditional Cayley trees (where terminals are of the same order) considered as an excitable media for input signal propagation has been recently used to investigate dynamic range \cite{Gollo2009,Gollo2013}. With our modelling approach for realistic neuron tree topology, it would be interesting to explore statistical properties of neuron trees in relation with dynamic ranges of neurons. In particular, it has been suggested that larger tree size could give larger dynamic ranges \cite{Gollo2009}; this can be tested with realistic virtual neuron trees using our modelling approach instead of traditional Cayley trees. Moreover, using more realistic virtual neuron trees also allow investigating the impact of tree asymmetry as well as tree shape in dynamic range.

\section*{Methods}
In this manuscript, we use 35 neuron digital reconstructions of pyramidal cells for analysis. Pyramidal neurons digital reconstructions are provided by authors in \cite{Zhou2018}. Briefly, neuron digita reconstructions were obtained as follows. The brain is from C57BL/6J mouse line, and its pyramid neurons are sparsely labelled with Adeno-associated virus (AAV). Fluorescence micro-optical sectioning tomography microscopy (fMOST) was used to image these labelled neurons which span different brain regions or even the whole brain. From the imaging dataset, GTree method was used to reconstruct these labelled neurons, all of which contained axonal and dendrite tree morphology. The animal experiments were approved by the Institutional Animal Ethics Committee of Huazhong University of Science and Technology, and all experiments were performed in accordance with relevant guidelines and regulations.


The digital reconstruction of each brain-wide neuron is stored in a morphological file (in ``SWC'' format) which includes position information of each traced point and its linkage points. For analysis, each empirically neuron is decomposed into one axonal tree and several dendritic trees; each tree stems from the same soma. The axon part is identified as the largest component (in terms of the number of branching points) among all subcomponents and the rest are considered as dendrites.

A subcomponent of digital reconstruction stemmed from the soma is excluded for analysis if there is no branching points in this subcomponent; this in fact corresponds to a tree of size $N=0$. Note that there is a small proportion of branching points which bifurcate to more than 2 branches; this could possibly due to errors in digital reconstruction. We manually adjust this by splitting such branching points so that all branching points bifurcate to 2 branches.


\begin{thebibliography}{10}
\expandafter\ifx\csname url\endcsname\relax
  \def\url#1{\texttt{#1}}\fi
\expandafter\ifx\csname urlprefix\endcsname\relax\def\urlprefix{URL }\fi
\expandafter\ifx\csname doiprefix\endcsname\relax\def\doiprefix{DOI }\fi
\providecommand{\bibinfo}[2]{#2}
\providecommand{\eprint}[2][]{\url{#2}}

\bibitem{Mainen1996}
\bibinfo{author}{Mainen, Z.~F.} \& \bibinfo{author}{Sejnowski, T.~J.}
\newblock \bibinfo{journal}{\bibinfo{title}{{Influence of dendritic structure
  on firing pattern in model neocortical neurons}}}.
\newblock {\emph{\JournalTitle{Nature}}} \textbf{\bibinfo{volume}{382}},
  \bibinfo{pages}{363--366} (\bibinfo{year}{1996}).

\bibitem{Vetter2001}
\bibinfo{author}{Vetter, P.}, \bibinfo{author}{Roth, A.} \&
  \bibinfo{author}{H{\"{a}}usser, M.}
\newblock \bibinfo{journal}{\bibinfo{title}{{Propagation of action potentials
  in dendrites depends on dendritic morphology.}}}
\newblock {\emph{\JournalTitle{Journal of neurophysiology}}}
  \textbf{\bibinfo{volume}{85}}, \bibinfo{pages}{926--937}
  (\bibinfo{year}{2001}).

\bibitem{Ferrante2013}
\bibinfo{author}{Ferrante, M.}, \bibinfo{author}{Migliore, M.} \&
  \bibinfo{author}{Ascoli, G.~A.}
\newblock \bibinfo{journal}{\bibinfo{title}{{Functional Impact of Dendritic
  Branch-Point Morphology}}}.
\newblock {\emph{\JournalTitle{Journal of Neuroscience}}}
  \textbf{\bibinfo{volume}{33}}, \bibinfo{pages}{2156--2165}
  (\bibinfo{year}{2013}).

\bibitem{Gollo2013}
\bibinfo{author}{Gollo, L.~L.}, \bibinfo{author}{Kinouchi, O.} \&
  \bibinfo{author}{Copelli, M.}
\newblock \bibinfo{journal}{\bibinfo{title}{{Single-neuron criticality
  optimizes analog dendritic computation}}}.
\newblock {\emph{\JournalTitle{Scientific Reports}}}
  \textbf{\bibinfo{volume}{3}}, \bibinfo{pages}{3222} (\bibinfo{year}{2013}).

\bibitem{Publio2012}
\bibinfo{author}{Publio, R.}, \bibinfo{author}{Ceballos, C.~C.} \&
  \bibinfo{author}{Roque, A.~C.}
\newblock \bibinfo{journal}{\bibinfo{title}{{Dynamic Range of Vertebrate Retina
  Ganglion Cells: Importance of Active Dendrites and Coupling by Electrical
  Synapses}}}.
\newblock {\emph{\JournalTitle{PLoS ONE}}} \textbf{\bibinfo{volume}{7}},
  \bibinfo{pages}{e48517} (\bibinfo{year}{2012}).

\bibitem{Yi2017}
\bibinfo{author}{Yi, G.~S.}, \bibinfo{author}{Wang, J.}, \bibinfo{author}{Deng,
  B.} \& \bibinfo{author}{Wei, X.~L.}
\newblock \bibinfo{journal}{\bibinfo{title}{{Morphology controls how
  hippocampal CA1 pyramidal neuron responds to uniform electric fields: A
  biophysical modeling study}}}.
\newblock {\emph{\JournalTitle{Scientific Reports}}}
  \textbf{\bibinfo{volume}{7}}, \bibinfo{pages}{3210} (\bibinfo{year}{2017}).

\bibitem{Ascoli2000}
\bibinfo{author}{Ascoli, G.~A.} \& \bibinfo{author}{Krichmar, J.~L.}
\newblock \bibinfo{journal}{\bibinfo{title}{{L-neuron: A modeling tool for the
  efficient generation and parsimonious description of dendritic morphology}}}.
\newblock {\emph{\JournalTitle{Neurocomputing}}}
  \textbf{\bibinfo{volume}{32-33}}, \bibinfo{pages}{1003--1011}
  (\bibinfo{year}{2000}).

\bibitem{Ascoli2001}
\bibinfo{author}{Ascoli, G.~A.}, \bibinfo{author}{Krichmar, J.~L.},
  \bibinfo{author}{Nasuto, S.~J.} \& \bibinfo{author}{Senft, S.~L.}
\newblock \bibinfo{journal}{\bibinfo{title}{{Generation, description and
  storage of dendritic morphology data}}}.
\newblock {\emph{\JournalTitle{Philosophical Transactions of the Royal Society
  B: Biological Sciences}}} \textbf{\bibinfo{volume}{356}},
  \bibinfo{pages}{1131--1145} (\bibinfo{year}{2001}).

\bibitem{Eberhard2006}
\bibinfo{author}{Eberhard, J.~P.}, \bibinfo{author}{Wanner, A.} \&
  \bibinfo{author}{Wittum, G.}
\newblock \bibinfo{journal}{\bibinfo{title}{{NeuGen: A tool for the generation
  of realistic morphology of cortical neurons and neural networks in 3D}}}.
\newblock {\emph{\JournalTitle{Neurocomputing}}} \textbf{\bibinfo{volume}{70}},
  \bibinfo{pages}{327--342} (\bibinfo{year}{2006}).

\bibitem{VanPelt1985}
\bibinfo{author}{{Van Pelt}, J.} \& \bibinfo{author}{Verwer, R.~W.}
\newblock \bibinfo{journal}{\bibinfo{title}{{Growth models (including terminal
  and segmental branching) for topological binary trees}}}.
\newblock {\emph{\JournalTitle{Bulletin of Mathematical Biology}}}
  \textbf{\bibinfo{volume}{47}}, \bibinfo{pages}{323--336}
  (\bibinfo{year}{1985}).

\bibitem{VanPelt1986a}
\bibinfo{author}{{Van Pelt}, J.} \& \bibinfo{author}{Verwer, R.~W.}
\newblock \bibinfo{journal}{\bibinfo{title}{{Topological properties of binary
  trees grown with order-dependent branching probabilities}}}.
\newblock {\emph{\JournalTitle{Bulletin of Mathematical Biology}}}
  \textbf{\bibinfo{volume}{48}}, \bibinfo{pages}{197--211}
  (\bibinfo{year}{1986}).

\bibitem{VanPelt2004}
\bibinfo{author}{{Van Pelt}, J.} \& \bibinfo{author}{Schierwagen, A.}
\newblock \bibinfo{journal}{\bibinfo{title}{{Morphological analysis and
  modeling of neuronal dendrites}}}.
\newblock {\emph{\JournalTitle{Mathematical Biosciences}}}
  \textbf{\bibinfo{volume}{188}}, \bibinfo{pages}{147--155}
  (\bibinfo{year}{2004}).

\bibitem{vanPelt2007}
\bibinfo{author}{{Van Pelt}, J.} \& \bibinfo{author}{Uylings, H. B.~M.}
\newblock \bibinfo{title}{{Modeling Neuronal Growth and Shape}}.
\newblock In \bibinfo{editor}{Laublicher, M.~D.} \&
  \bibinfo{editor}{M{\"{u}}ller, G.~B.} (eds.)
  \emph{\bibinfo{booktitle}{Modeling Biology – Structures, Behaviors,
  Evolution}}, \bibinfo{pages}{195--215} (\bibinfo{publisher}{The MIT Press},
  \bibinfo{address}{Cambridge, Massachusetts}, \bibinfo{year}{2007}).

\bibitem{Koene2009}
\bibinfo{author}{Koene, R.~A.} \emph{et~al.}
\newblock \bibinfo{journal}{\bibinfo{title}{{NETMORPH: A framework for the
  stochastic generation of large scale neuronal networks with realistic neuron
  morphologies}}}.
\newblock {\emph{\JournalTitle{Neuroinformatics}}}
  \textbf{\bibinfo{volume}{7}}, \bibinfo{pages}{195--210}
  (\bibinfo{year}{2009}).

\bibitem{Cuntz2010}
\bibinfo{author}{Cuntz, H.}, \bibinfo{author}{Forstner, F.},
  \bibinfo{author}{Borst, A.} \& \bibinfo{author}{H{\"{a}}usser, M.}
\newblock \bibinfo{journal}{\bibinfo{title}{{One rule to grow them all: A
  general theory of neuronal branching and its practical application}}}.
\newblock {\emph{\JournalTitle{PLoS Computational Biology}}}
  \textbf{\bibinfo{volume}{6}}, \bibinfo{pages}{e1000877}
  (\bibinfo{year}{2010}).

\bibitem{Cuntz2012}
\bibinfo{author}{Cuntz, H.}, \bibinfo{author}{Mathy, A.} \&
  \bibinfo{author}{Hausser, M.}
\newblock \bibinfo{journal}{\bibinfo{title}{{A scaling law derived from optimal
  dendritic wiring}}}.
\newblock {\emph{\JournalTitle{Proceedings of the National Academy of
  Sciences}}} \textbf{\bibinfo{volume}{109}}, \bibinfo{pages}{11014--11018}
  (\bibinfo{year}{2012}).

\bibitem{Luczak2010}
\bibinfo{author}{Luczak, A.}
\newblock \bibinfo{journal}{\bibinfo{title}{{Measuring neuronal branching
  patterns using model-based approach}}}.
\newblock {\emph{\JournalTitle{Frontiers in Computational Neuroscience}}}
  \textbf{\bibinfo{volume}{4}}, \bibinfo{pages}{135} (\bibinfo{year}{2010}).

\bibitem{Luczak2006}
\bibinfo{author}{Luczak, A.}
\newblock \bibinfo{journal}{\bibinfo{title}{{Spatial embedding of neuronal
  trees modeled by diffusive growth}}}.
\newblock {\emph{\JournalTitle{Journal of Neuroscience Methods}}}
  \textbf{\bibinfo{volume}{157}}, \bibinfo{pages}{132--141}
  (\bibinfo{year}{2006}).

\bibitem{Gillette2015b}
\bibinfo{author}{Gillette, T.~A.} \& \bibinfo{author}{Ascoli, G.~A.}
\newblock \bibinfo{journal}{\bibinfo{title}{{Topological characterization of
  neuronal arbor morphology via sequence representation: II - global
  alignment}}}.
\newblock {\emph{\JournalTitle{BMC Bioinformatics}}}
  \textbf{\bibinfo{volume}{16}}, \bibinfo{pages}{209} (\bibinfo{year}{2015}).

\bibitem{Mohan2015}
\bibinfo{author}{Mohan, H.} \emph{et~al.}
\newblock \bibinfo{journal}{\bibinfo{title}{{Dendritic and Axonal Architecture
  of Individual Pyramidal Neurons across Layers of Adult Human Neocortex}}}.
\newblock {\emph{\JournalTitle{Cerebral Cortex}}}
  \textbf{\bibinfo{volume}{25}}, \bibinfo{pages}{4839--4853}
  (\bibinfo{year}{2015}).

\bibitem{Li2010}
\bibinfo{author}{Li, A.} \emph{et~al.}
\newblock \bibinfo{journal}{\bibinfo{title}{{Micro-optical sectioning
  tomography to obtain a high-resolution atlas of the mouse brain.}}}
\newblock {\emph{\JournalTitle{Science}}} \textbf{\bibinfo{volume}{330}},
  \bibinfo{pages}{1404--1408} (\bibinfo{year}{2010}).

\bibitem{Ragan2012}
\bibinfo{author}{Ragan, T.} \emph{et~al.}
\newblock \bibinfo{journal}{\bibinfo{title}{{Serial two-photon tomography for
  automated ex vivo mouse brain imaging}}}.
\newblock {\emph{\JournalTitle{Nature Methods}}} \textbf{\bibinfo{volume}{9}},
  \bibinfo{pages}{255--258} (\bibinfo{year}{2012}).

\bibitem{Jefferis2012}
\bibinfo{author}{Jefferis, G.~S.} \& \bibinfo{author}{Livet, J.}
\newblock \bibinfo{journal}{\bibinfo{title}{{Sparse and combinatorial neuron
  labelling}}}.
\newblock {\emph{\JournalTitle{Current Opinion in Neurobiology}}}
  \textbf{\bibinfo{volume}{22}}, \bibinfo{pages}{101--110}
  (\bibinfo{year}{2012}).

\bibitem{Chung2013}
\bibinfo{author}{Chung, K.} \& \bibinfo{author}{Deisseroth, K.}
\newblock \bibinfo{journal}{\bibinfo{title}{{CLARITY for mapping the nervous
  system}}}.
\newblock {\emph{\JournalTitle{Nature Methods}}} \textbf{\bibinfo{volume}{10}},
  \bibinfo{pages}{508--513} (\bibinfo{year}{2013}).

\bibitem{Zhou2018}
\bibinfo{author}{Zhou, H.} \emph{et~al.}
\newblock \bibinfo{journal}{\bibinfo{title}{{Dense reconstruction of brain-wide
  neuronal population close to the ground truth}}}.
\newblock {\emph{\JournalTitle{bioRxiv}}} \bibinfo{pages}{223834}
  (\bibinfo{year}{2018}).

\bibitem{Quan2016}
\bibinfo{author}{Quan, T.} \emph{et~al.}
\newblock \bibinfo{journal}{\bibinfo{title}{{NeuroGPS-Tree: automatic
  reconstruction of large-scale neuronal populations with dense neurites}}}.
\newblock {\emph{\JournalTitle{Nature Methods}}} \textbf{\bibinfo{volume}{13}},
  \bibinfo{pages}{51--54} (\bibinfo{year}{2016}).

\bibitem{Vormberg2017}
\bibinfo{author}{Vormberg, A.}, \bibinfo{author}{Effenberger, F.},
  \bibinfo{author}{Muellerleile, J.} \& \bibinfo{author}{Cuntz, H.}
\newblock \bibinfo{journal}{\bibinfo{title}{{Universal features of dendrites
  through centripetal branch ordering}}}.
\newblock {\emph{\JournalTitle{PLoS Computational Biology}}}
  \textbf{\bibinfo{volume}{13}}, \bibinfo{pages}{e1005615}
  (\bibinfo{year}{2017}).

\bibitem{Gillette2015a}
\bibinfo{author}{Gillette, T.~A.} \& \bibinfo{author}{Ascoli, G.~A.}
\newblock \bibinfo{journal}{\bibinfo{title}{{Topological characterization of
  neuronal arbor morphology via sequence representation: I - motif analysis}}}.
\newblock {\emph{\JournalTitle{BMC Bioinformatics}}}
  \textbf{\bibinfo{volume}{16}}, \bibinfo{pages}{216} (\bibinfo{year}{2015}).

\bibitem{Molnar2006}
\bibinfo{author}{Molnar, P.}
\newblock \bibinfo{journal}{\bibinfo{title}{{On geometrical scaling of Cayley
  trees and river networks}}}.
\newblock {\emph{\JournalTitle{Journal of Hydrology}}}
  \textbf{\bibinfo{volume}{322}}, \bibinfo{pages}{199--210}
  (\bibinfo{year}{2006}).

\bibitem{Sholl1953}
\bibinfo{author}{Sholl, D.~A.}
\newblock \bibinfo{journal}{\bibinfo{title}{{Dendritic organization in the
  neurons of the visual and motor cortices of the cat}}}.
\newblock {\emph{\JournalTitle{Journal of Anatomy}}}
  \textbf{\bibinfo{volume}{87}}, \bibinfo{pages}{387--406}
  (\bibinfo{year}{1953}).

\bibitem{VanPelt1992}
\bibinfo{author}{{Van Pelt}, J.}, \bibinfo{author}{Uylings, H.~B.},
  \bibinfo{author}{Verwer, R.~W.}, \bibinfo{author}{Pentney, R.~J.} \&
  \bibinfo{author}{Woldenberg, M.~J.}
\newblock \bibinfo{journal}{\bibinfo{title}{{Tree asymmetry-A sensitive and
  practical measure for binary topological trees}}}.
\newblock {\emph{\JournalTitle{Bulletin of Mathematical Biology}}}
  \textbf{\bibinfo{volume}{54}}, \bibinfo{pages}{759--784}
  (\bibinfo{year}{1992}).

\bibitem{Samsonovich2006}
\bibinfo{author}{Samsonovich, A.~V.} \& \bibinfo{author}{Ascoli, G.~A.}
\newblock \bibinfo{journal}{\bibinfo{title}{{Morphological homeostasis in
  cortical dendrites.}}}
\newblock {\emph{\JournalTitle{Proceedings of the National Academy of Sciences
  of the United States of America}}} \textbf{\bibinfo{volume}{103}},
  \bibinfo{pages}{1569--1574} (\bibinfo{year}{2006}).

\bibitem{Samsonovich2005}
\bibinfo{author}{Samsonovich, A.~V.} \& \bibinfo{author}{Ascoli, G.~A.}
\newblock \bibinfo{journal}{\bibinfo{title}{{Statistical determinants of
  dendritic morphology in hippocampal pyramidal neurons: A hidden Markov
  model}}}.
\newblock {\emph{\JournalTitle{Hippocampus}}} \textbf{\bibinfo{volume}{15}},
  \bibinfo{pages}{166--183} (\bibinfo{year}{2005}).

\bibitem{Ostilli2012}
\bibinfo{author}{Ostilli, M.}
\newblock \bibinfo{journal}{\bibinfo{title}{{Cayley Trees and Bethe Lattices: A
  concise analysis for mathematicians and physicists}}}.
\newblock {\emph{\JournalTitle{Physica A: Statistical Mechanics and its
  Applications}}} \textbf{\bibinfo{volume}{391}}, \bibinfo{pages}{3417--3423}
  (\bibinfo{year}{2012}).

\bibitem{Albert2002}
\bibinfo{author}{Albert, R.} \& \bibinfo{author}{Barab{\'{a}}si, A.-L.}
\newblock \bibinfo{journal}{\bibinfo{title}{{Statistical mechanics of complex
  networks}}}.
\newblock {\emph{\JournalTitle{Reviews of Modern Physics}}}
  \textbf{\bibinfo{volume}{74}}, \bibinfo{pages}{47--97}
  (\bibinfo{year}{2002}).

\bibitem{Donohue2008}
\bibinfo{author}{Donohue, D.~E.} \& \bibinfo{author}{Ascoli, G.~A.}
\newblock \bibinfo{journal}{\bibinfo{title}{{A comparative computer simulation
  of dendritic morphology}}}.
\newblock {\emph{\JournalTitle{PLoS Computational Biology}}}
  \textbf{\bibinfo{volume}{4}}, \bibinfo{pages}{e1000089}
  (\bibinfo{year}{2008}).

\bibitem{Zeng2017}
\bibinfo{author}{Zeng, H.} \& \bibinfo{author}{Sanes, J.~R.}
\newblock \bibinfo{title}{{Neuronal cell-type classification: Challenges,
  opportunities and the path forward}} (\bibinfo{year}{2017}).

\bibitem{Joo2013}
\bibinfo{author}{Joo, W.~J.}, \bibinfo{author}{Sweeney, L.~B.},
  \bibinfo{author}{Liang, L.} \& \bibinfo{author}{Luo, L.}
\newblock \bibinfo{journal}{\bibinfo{title}{{Linking cell fate, trajectory
  choice, and target selection: Genetic analysis of sema-2b in olfactory axon
  targeting}}}.
\newblock {\emph{\JournalTitle{Neuron}}} \textbf{\bibinfo{volume}{78}},
  \bibinfo{pages}{673--686} (\bibinfo{year}{2013}).

\bibitem{Han2018}
\bibinfo{author}{Han, Y.} \emph{et~al.}
\newblock \bibinfo{journal}{\bibinfo{title}{{The logic of single-cell
  projections from visual cortex}}}.
\newblock {\emph{\JournalTitle{Nature}}} \textbf{\bibinfo{volume}{556}},
  \bibinfo{pages}{51--56} (\bibinfo{year}{2018}).

\bibitem{Schubert2006}
\bibinfo{author}{Schubert, D.}, \bibinfo{author}{K{\"{o}}tter, R.},
  \bibinfo{author}{Luhmann, H.~J.} \& \bibinfo{author}{Staiger, J.~F.}
\newblock \bibinfo{journal}{\bibinfo{title}{{Morphology, electrophysiology and
  functional input connectivity of pyramidal neurons characterizes a genuine
  layer Va in the primary somatosensory cortex}}}.
\newblock {\emph{\JournalTitle{Cerebral Cortex}}}
  \textbf{\bibinfo{volume}{16}}, \bibinfo{pages}{223--236}
  (\bibinfo{year}{2006}).

\bibitem{Marx2013}
\bibinfo{author}{Marx, M.} \& \bibinfo{author}{Feldmeyer, D.}
\newblock \bibinfo{journal}{\bibinfo{title}{{Morphology and physiology of
  excitatory neurons in layer 6b of the somatosensory rat barrel cortex.}}}
\newblock {\emph{\JournalTitle{Cerebral cortex}}}
  \textbf{\bibinfo{volume}{23}}, \bibinfo{pages}{2803--2817}
  (\bibinfo{year}{2013}).

\bibitem{Gollo2009}
\bibinfo{author}{Gollo, L.~L.}, \bibinfo{author}{Kinouchi, O.} \&
  \bibinfo{author}{Copelli, M.}
\newblock \bibinfo{journal}{\bibinfo{title}{{Active dendrites enhance neuronal
  dynamic range}}}.
\newblock {\emph{\JournalTitle{PLoS Computational Biology}}}
  \textbf{\bibinfo{volume}{5}}, \bibinfo{pages}{e1000402}
  (\bibinfo{year}{2009}).

\end{thebibliography}

\section*{Acknowledgements}
CL acknowledges financial support from National foundation of Science in China (NFSC, grant No. 11701201 and No. 11871061), YZ acknowledges financial support from by NFSC (grant No. 11701200 and 11871262) and TQ acknowledges financial support from Science Fund for Young and Middle-aged Creative Research Group of the Universities in Hubei Province (Grant No. T201520) and NSFC (grant No. 81771913). We also thank the Optical Bioimaging Core Facility of WNLO-HUST for the support in data acquisition, and the Analytical and Testing Center of HUST for spectral measurements.  We thank Prof. Shaoqun Zeng for helpful discussions on actual neuron data.

\section*{Author contributions statement}
YH, CL and TQ performed the analysis, CL and YH performed modelling, and YZ and CL conceived and designed the project. All authors reviewed the manuscript.

\section*{Additional information}
\textbf{Competing interests}: The authors declare no competing interests.

\section*{Data availability}
The data that support the findings of this study are available from the corresponding author on request.

\end{document}